\documentclass[reqno]{article}%
\usepackage{amsmath}
\usepackage{graphicx}
\usepackage{amsfonts}
\usepackage{amssymb}%
\newtheorem{theorem}{Theorem}

\newtheorem{claim}[theorem]{Claim}

\newtheorem{corollary}[theorem]{Corollary}

\newtheorem{lemma}[theorem]{Lemma}

\newtheorem{proposition}[theorem]{Proposition}

\begin{document}

\title{Book Ramsey numbers I}
\author{Vladimir Nikiforov\thanks{Department of Mathematical Sciences, University of
Memphis, Memphis TN 38152, USA} \ and Cecil Rousseau$^{\ast}$}
\maketitle

\begin{abstract}
A book $B_{p}$ is a graph consisting of $p$ triangles sharing a common edge.
In this paper we prove that if $p\leq q/6-o\left(  q\right)  $ and $q$ is
large then the Ramsey number $r\left(  B_{p},B_{q}\right)  $ is given by%
\[
r\left(  B_{p},B_{q}\right)  =2q+3
\]
and the constant $1/6$ is essentially best possible.

Our proof is based on Szemer\'{e}di's uniformity lemma and a stability result
for books.

\textbf{Keywords: }\textit{book, Ramsey number, uniformity lemma}

\end{abstract}

\section{Introduction}

Our notation and terminology are standard (see, e.g., \cite{Bol}). Thus,
$G\left(  n,m\right)  $ is a graph of order $n$ and size $m;$ for a graph $G$
and a vertex $u\in V\left(  G\right)  $ we write $\Gamma_{G}\left(  u\right)
$ for the set of vertices adjacent to $u;$ $d_{G}\left(  u\right)  =\left|
\Gamma\left(  u\right)  \right|  $ is the degree of $u;$ we write $d\left(
u\right)  $ and $\Gamma\left(  u\right)  $ instead of $d_{G}\left(  u\right)
$ and $\Gamma_{G}\left(  u\right)  $ when it is clear which graph $G$ is to be taken.

Unless explicitely stated, all graphs are assumed to be defined on the vertex
set $\left[  n\right]  =\left\{  1,2,...,n\right\}  .$

A \emph{book }of size $q$ consists of $q$ triangles sharing a common edge. We
write $bk\left(  G\right)  $ for the size of the largest book in a graph $G$
and call it the \emph{booksize }of\emph{ }$G.$

Books have attracted considerable attention in Ramsey graph theory (see, e.g.,
\cite{RS}, \cite{FRS82} and \cite{NR}).

The study of Ramsey numbers for books was initiated by Rousseau and Sheehan in
\cite{RS} where the following was proved.

\begin{theorem}
[Rousseau and Sheehan]For all $q>1$,
\[
r(B_{1},B_{q})=2q+3.
\]

\end{theorem}

Somewhat later Faudree, Sheehan and Rousseau strengthened this in \cite{FRS82}
in the following way.

\begin{theorem}
[Faudree et al.]Suppose $1\leq p\leq q.$ Then%
\[
r(B_{p},B_{q})=2q+3.
\]
for all%
\[
q\geq\left(  p-1\right)  \left(  16p^{3}+16p^{2}-24p-10\right)  +1.
\]

\end{theorem}

The quartic constraint of $q$ on $p$ was in turn reduced to linear by
Nikiforov and Rousseau in \cite{NR}.

\begin{theorem}
[Nikiforov, Rousseau]There exists a positive constant $c$ such that
\[
r(B_{p},B_{q})=2q+3
\]
for all $q\geq cp$.
\end{theorem}

In \cite{NR} it was found that $c\geq10^{-6}$. In fact it turns out that
\[
c=1/6+o\left(  1\right)
\]
and the proof of this inequality is our main goal in this chapter.

\begin{theorem}
\label{th1} For every $\varepsilon$ $\left(  0<\varepsilon<1/6\right)  $ if
$q$ is large and $p\leq\left(  1/6-\varepsilon\right)  q$ then
\[
r\left(  B_{p},B_{q}\right)  =2q+3.
\]

\end{theorem}

Taking the union of two disjoint complete graphs of order $q+1$ we immediately
see that
\[
r\left(  B_{p},B_{q}\right)  \geq2q+3,
\]
so all we have to prove is that, under the assumptions of the theorem, the
inequality
\begin{equation}
r\left(  B_{p},B_{q}\right)  \leq2q+3 \label{mainin}%
\end{equation}
holds. The proof is obtained essentially as a corollary of theorem that is
considered in the following section.

\section{A Ramsey type stability result}

We shall prove the following.

\begin{theorem}
\label{th2} There exists a constant $c>0$ such that for every $\xi$ with
$0<\xi<c$ and for every graph $G$ of sufficiently large order $n,$ one of the
following holds:

\emph{(i)} $bk\left(  \overline{G}\right)  >n/2;$

\emph{(ii)} $bk\left(  G\right)  >\left(  1/12-10^{-6}\xi^{6}\right)  n;$

\emph{(iii)} $G$ contains an induced bipartite graph $G_{0}$ of order at least
$\left(  1-\xi\right)  n$ and
\[
\delta\left(  G_{0}\right)  >\left(  \frac{1}{2}-2\xi\right)  n.
\]

\end{theorem}

First we shall state some preliminary results related to Szemer\'{e}di's
uniformity lemma and a stability theorem for books in graphs.

\subsection{Results related to Szemer\'{e}di's uniformity lemma}

For notation and definitions related to Szemer\'{e}di's uniformity lemma see,
e.g. \cite{KoSi}. We also need a few technical results; the first one is a
basic property of $\varepsilon$-uniform pairs (see \cite{KoSi}, Fact 1.4).

The proofs of the stated results are collected in section \ref{secProofs}.

\begin{lemma}
\label{XPle} Suppose $0<\varepsilon<d\leq1$ and $\left(  A,B\right)  $ is an
$\varepsilon$-uniform pair with $e\left(  A,B\right)  =d\left|  A\right|
\left|  B\right|  .$ Then there are at most $2\varepsilon\left|  A\right|
^{2}$ sets $\left\{  u,v\right\}  \subset A$ with%
\[
\left|  \Gamma\left(  u\right)  \cap\Gamma\left(  v\right)  \cap B\right|
\leq\left(  d-\varepsilon\right)  ^{2}\left|  B\right|  .
\]
$\ \hfill\square$
\end{lemma}

The next lemma gives a lower bound on the number of triangles in a graph that
consists of several $\varepsilon$-uniform pairs sharing a common part.

\begin{lemma}
\label{dle} Let $\varepsilon>0$ and $H$ be a graph whose vertices are
partitioned as
\[
V\left(  H\right)  =A\cup B_{1}\cup...\cup B_{k}%
\]
with
\[
\left|  A\right|  =\left|  B_{1}\right|  =...=\left|  B_{k}\right|  =t.
\]
Suppose that for every $i\in\left[  k\right]  $ the pair $\left(
A,B_{i}\right)  $ is $\varepsilon$-uniform and $e\left(  A,B_{i}\right)  \geq
d_{i}t^{2}$. Then there are at least
\[
t\left(  e\left(  A\right)  -2\varepsilon t^{2}\right)  \sum_{i=1}^{k}%
d_{i}^{2}-2\varepsilon kte\left(  A\right)
\]
triangles in $H$ having exactly $2$ vertices in $A.$
\end{lemma}

By averaging over the edges in $A$ we obtain the following corollary.

\begin{corollary}
\label{bs1} Under the conditions of Lemma \ref{dle}, if $e\left(  A\right)
>0$ then there is a book whose base is in $A$ and whose size is at least%
\[
t\left(  1-\frac{2\varepsilon t^{2}}{e\left(  A\right)  }\right)  \sum
_{i=1}^{k}d_{i}^{2}-2\varepsilon kt.
\]

\end{corollary}

Similar results hold for books whose bases belong to different blobs.

\begin{lemma}
\label{XPle2} Suppose $\varepsilon,$ $d_{1},$ $d_{2}$ are real numbers with
\[
0<2\varepsilon\leq d_{1}\leq1,\text{ \ \ \ \ }2\varepsilon\leq d_{2}\leq1.
\]
Let $\left(  A_{1},B\right)  ,$ $\left(  A_{2},B\right)  $ be $\varepsilon
$-uniform pairs with
\[
e\left(  A_{i},B\right)  =d_{i}\left|  A_{i}\right|  \left|  B\right|  ,\text{
\ \ }\left(  i=1,2\right)  .
\]
Then there are at most $2\varepsilon\left|  A_{1}\right|  \left|
A_{2}\right|  $ pairs $\left(  u,v\right)  $ with $u\in A_{1},$ $v\in A_{2}$
such that
\begin{equation}
\left|  \Gamma\left(  u\right)  \cap\Gamma\left(  v\right)  \cap B\right|
\leq\left(  d_{1}-\varepsilon\right)  \left(  d_{2}-\varepsilon\right)
\left|  B\right|  . \label{i03}%
\end{equation}

\end{lemma}

\begin{lemma}
\label{lebs} Suppose $\varepsilon$ $>0$ and $H$ is a graph whose vertices are
partitioned as
\[
V\left(  H\right)  =A_{1}\cup A_{2}\cup B_{1}\cup...\cup B_{k}%
\]
with
\[
\left|  A_{1}\right|  =\left|  A_{2}\right|  =\left|  B_{1}\right|
=...=\left|  B_{k}\right|  =t.
\]
Suppose that for every $i\in\left[  2\right]  ,$ $j\in\left[  k\right]  $ the
pairs $\left(  A_{i},B_{j}\right)  $ are $\varepsilon$-uniform and
\[
e\left(  A_{i},B_{j}\right)  \geq d_{ij}t^{2},
\]
Then there are at least
\[
t\left(  e\left(  A_{1},A_{2}\right)  -2\varepsilon t^{2}\right)  \left(
\sum_{i=1}^{k}d_{1i}d_{2i}\right)  -2\varepsilon kte\left(  A_{1}%
,A_{2}\right)
\]
triangles in $H$ having exactly one vertex in $A_{1}$ and one vertex in
$A_{2}.$
\end{lemma}

By averaging over the edges in $E\left(  A_{1},A_{2}\right)  $ we obtain the
following corollary.

\begin{corollary}
\label{bs2} Under the conditions of Lemma \ref{lebs}, if $e\left(  A_{1}%
,A_{2}\right)  >0$ then there is a book whose base belongs to $E\left(
A_{1},A_{2}\right)  $ and whose size is at least
\[
t\left(  1-\frac{2\varepsilon t^{2}}{e\left(  A_{1},A_{2}\right)  }\right)
\left(  \sum_{i=1}^{k}d_{1i}d_{2i}\right)  -2\varepsilon kt.
\]

\end{corollary}

\subsection{A stability theorem for books}

A key ingredient in our proof of Theorem \ref{th2} will be the following
result that was proved in \cite{BoNi}.

\begin{theorem}
\label{stab} For every $\alpha$ with $0<\alpha<10^{-5}$ and every graph
$G=G\left(  n\right)  $ with
\begin{equation}
e\left(  G\right)  \geq\left(  \frac{1}{4}-\alpha\right)  n^{2} \label{cond1}%
\end{equation}
either
\[
bk\left(  G\right)  >\left(  \frac{1}{6}-2\alpha^{1/3}\right)  n
\]
or $G$ contains an induced bipartite graph $G_{1}$ of order at least $\left(
1-\alpha^{1/3}\right)  n$ and with minimal degree
\[
\delta\left(  G_{1}\right)  \geq\left(  \frac{1}{2}-4\alpha^{1/3}\right)  n.
\]
$\ \hfill\square$
\end{theorem}

One can see immediately that this theorem has a close relationship to Theorem
\ref{th2}.

\subsection{Proof of Theorem \ref{th2}}

Instead of the graph $G$ and its complement we shall consider a blue-red
coloring of $K_{n}$ so that the blue edges correspond to the edges of $G$.

Set
\[
bk_{R}=bk\left(  \overline{G}\right)  \text{ \ \ \ \ \ }bk_{B}=bk\left(
G\right)  .
\]

If $X\subset\left[  n\right]  $ then $e_{R}\left(  X\right)  $ and
$e_{B}\left(  X\right)  $ denote respectively the number of the red and blue
edges induced by $X.$ Similarly if $X,Y\subset\left[  n\right]  $ are disjoint
sets then $e_{R}\left(  X,Y\right)  $ and $e_{B}\left(  X,Y\right)  $ denote
repectively the number of the red and blue $X-Y$ edges.

Assume for every choice of $\xi>0$ we have for large enough $n,$%
\begin{equation}
bk_{R}\leq\frac{n}{2},\text{ \ \ \ \ }bk_{B}\leq\left(  \frac{1}{12}%
-10^{-6}\xi^{6}\right)  n. \label{assum0}%
\end{equation}

Our goal is to show that these conditions imply (iii). To achieve this we
shall apply Theorem \ref{stab}, but to do so we have to ensure that the number
of blue edges is sufficiently close to $n^{2}/4,$ so that assumption
(\ref{cond1}) holds. The bulk of our proof is dedicated to this purpose. We
shall use Szemer\'{e}di's uniformity lemma to derive a number of conditions on
our edge coloring implying eventually that the number of blue edges is
sufficiently close to $n^{2}/4$ and then by Theorem \ref{stab} we shall
complete the proof.

Fix a small enough value $\xi>0$ and set
\[
\beta=\frac{1}{17}\left(  \frac{\xi}{2}\right)  ^{3},\text{ \ \ \ }%
\gamma=\beta^{3},\text{ \ \ \ }\varepsilon=\frac{\gamma^{2}}{2}.
\]

In the course of our proof we shall frequently use the fact that
\[
\varepsilon\ll\gamma\ll\beta\ll\xi,
\]
so, selecting $\xi$ sufficiently small, we can make $\varepsilon,$ $\gamma$
and $\beta$ as small as we like.

According to the uniformity lemma of Szemer\'{e}di for $n$ sufficiently large
there exists a partition
\[
\left[  n\right]  =V_{0}\cup V_{1}\cup...\cup V_{k},
\]
so that
\[
\left|  V_{0}\right|  <\varepsilon n,\text{ \ \ \ }\left|  V_{1}\right|
=...=\left|  V_{k}\right|
\]
and all but $\varepsilon k^{2}$ pairs $\left(  V_{i},V_{j}\right)  $ are
$\varepsilon$-uniform. As usual, we shall call the sets $V_{1},...,V_{k}$
\emph{blobs}.

In addition we may and shall suppose that $k$ is sufficiently large and for
every $i\in\left[  k\right]  ,$ less than $\varepsilon k$ pairs $\left(
V_{i},V_{j}\right)  $ are not $\varepsilon$-uniform.

Set $\left|  V_{1}\right|  =t;$ in view of $17^{-3}2^{-6}>10^{-6},$ the
assumption (\ref{assum0}) implies%
\begin{equation}
bk_{R}\leq\frac{1}{2}\left(  1+\varepsilon\right)  kt,\text{ \ \ \ \ }%
bk_{B}\leq\left(  \frac{1}{12}-\gamma\right)  \left(  1+\varepsilon\right)
kt. \label{assum}%
\end{equation}

For every $i,j\in\left[  k\right]  ,$ $\left(  i\neq j\right)  ,$ let $d_{ij}$
be the red density of the pair $\left(  V_{i},V_{j}\right)  ,$ i.e.,%
\[
d_{ij}=e_{R}\left(  V_{i},V_{j}\right)  /t^{2}.
\]

We shall prove four claims showing that (\ref{assum}) imposes rigid structural
restrictions on our edge coloring. The proofs of these claims are
straightforward but technical, so to keep the main line clear we have
collected them in section \ref{secProofs}.

First we shall prove that no blob contains significantly many edges of both colors.

\begin{claim}
\label{cl1} For every $i,$ either $e_{R}\left(  V_{i}\right)  <\gamma t^{2}$
or $e_{B}\left(  V_{i}\right)  <\gamma t^{2}.$
\end{claim}

We call a blob $\emph{blue}$ if it induces at most $\gamma t^{2}$ red edges,
and \emph{red }if it induces at most $\gamma t^{2}$ blue edges. Observe that
by Claim \ref{cl1} every blob is either red or blue but, of course, none is
both red and blue.

Next we shall prove that there are no three blobs $V_{i},$ $V_{j},$ $V_{l}$
all having significantly many red edges and each two joined by significantly
many blue edges.

\begin{claim}
\label{cl2}\textbf{ }There are no three blobs $V_{i},$ $V_{j},$ $V_{l}$\emph{
}such that%
\[
e_{R}\left(  V_{i}\right)  \geq\gamma t^{2},\text{ }e_{R}\left(  V_{j}\right)
\geq\gamma t^{2},\text{ }e_{R}\left(  V_{l}\right)  \geq\gamma t^{2}%
\]
and
\[
e_{B}\left(  V_{i},V_{j}\right)  \geq\gamma t^{2},\text{ }e_{B}\left(
V_{i},V_{l}\right)  \geq\gamma t^{2},\text{ }e_{B}\left(  V_{l},V_{j}\right)
\geq\gamma t^{2}.
\]

\end{claim}

The next claim shows that there are no two blobs containing significantly many
blue edges joined by significantly many red edges.

\begin{claim}
\label{cl3} There are no two blobs $V_{i},$ $V_{j}$\emph{ }such that%
\[
e_{B}\left(  V_{i}\right)  \geq\gamma t^{2},\text{ }e_{B}\left(  V_{j}\right)
\geq\gamma t^{2}%
\]
and
\[
e_{R}\left(  V_{i},V_{j}\right)  \geq\gamma t^{2}.
\]

\end{claim}

The next claim shows that all blobs are red.

\begin{claim}
\label{cl4} There are no blue blobs
\end{claim}

To finish the proof we shall show that the number of blue edges is arbitrarily
close to $n^{2}/4.$

Recall that $d_{ij}$ is the red density of the pair $\left(  V_{i}%
,V_{j}\right)  .$ Define the graphs $H_{irr},$ $H_{blue},$ $H_{mid}$ and
$H_{red}$ on the vertex set $\left[  k\right]  $ as follows:

(a) $\left(  i,j\right)  \in E\left(  H_{irr}\right)  $ iff the pair $\left(
V_{i},V_{j}\right)  $ is not $\varepsilon$-uniform;

(b) $\left(  i,j\right)  \in E\left(  H_{blue}\right)  $ iff the pair $\left(
V_{i},V_{j}\right)  $ is $\varepsilon$-uniform and
\[
d_{ij}<\beta;
\]

(c) $\left(  i,j\right)  \in E\left(  H_{mid}\right)  $ iff the pair $\left(
V_{i},V_{j}\right)  $ is $\varepsilon$-uniform and
\[
\beta\leq d_{ij}<1-\gamma;
\]

(d) $\left(  i,j\right)  \in E\left(  H_{red}\right)  $ iff the pair $\left(
V_{i},V_{j}\right)  $ is $\varepsilon$-uniform and
\[
d_{ij}\geq1-\gamma;
\]
Observe that the graphs $H_{irr},$ $H_{blue},$ $H_{mid}$ and $H_{red}$ are
pairwise edge disjoint.

Let $i$ be any vertex in $H_{red}.$ Estimating the average size of the red
books whose base is in $E\left(  V_{i}\right)  $ we obtain%
\begin{align}
\left(  d_{H_{red}}\left(  i\right)  -\varepsilon k\right)  \left(
1-\varepsilon-\gamma\right)  ^{2}t+\left(  d_{H_{mid}}\left(  i\right)
-\varepsilon k\right)  \left(  \beta-\varepsilon\right)  ^{2}  &  \leq
bk_{R}\label{in1}\\
&  \leq\frac{1}{2}\left(  1+\varepsilon\right)  kt\nonumber
\end{align}
and hence%
\begin{equation}
\left(  \frac{2e\left(  H_{red}\right)  }{k^{2}}-\varepsilon\right)  \left(
1-\varepsilon-\gamma\right)  ^{2}+\left(  \frac{2e\left(  H_{mid}\right)
}{k^{2}}-\varepsilon\right)  \left(  \beta-\varepsilon\right)  ^{2}\leq
\frac{1}{2}+\gamma. \label{in2}%
\end{equation}
Since by Claim \ref{cl2} the complement of $H_{red}$ is triangle-free, by
Tur\'{a}n's theorem we have%
\[
e\left(  H_{red}\right)  \geq\binom{k}{2}-\frac{k^{2}}{4}>\left(  \frac{1}%
{4}-\varepsilon\right)  k^{2}%
\]
for $k$ sufficiently large. Hence, we see that
\begin{align*}
\left(  \frac{2e\left(  H_{red}\right)  }{k^{2}}-\varepsilon\right)  \left(
1-\varepsilon-\gamma\right)  ^{2}  &  \geq\left(  \frac{1}{2}-3\varepsilon
\right)  \left(  1-\varepsilon-\gamma\right)  ^{2}\\
&  >\left(  \frac{1}{2}-3\gamma\right)  \left(  1-4\gamma\right)  >\frac{1}%
{2}-5\gamma.
\end{align*}
Therefore, from (\ref{in2}) we find that
\[
\left(  \frac{2e\left(  H_{mid}\right)  }{k^{2}}-\varepsilon\right)
\frac{\beta^{2}}{4}\leq\left(  \frac{2e\left(  H_{mid}\right)  }{k^{2}%
}-\varepsilon\right)  \left(  \beta-\varepsilon\right)  ^{2}\leq6\gamma
=6\beta^{3},
\]
and thus,%
\[
e\left(  H_{mid}\right)  <\left(  12\beta+\frac{\varepsilon}{2}\right)
k^{2}<13\beta k^{2}.
\]
On the other hand, from (\ref{in1}), we immediately have%
\[
\left(  d_{H_{red}}\left(  i\right)  -\varepsilon k\right)  \left(
1-\varepsilon-\gamma\right)  ^{2}\leq\frac{1}{2}\left(  1+\varepsilon\right)
k,
\]
and therefore,
\[
e\left(  H_{red}\right)  \leq\frac{1}{\left(  1-\varepsilon-\gamma\right)
^{2}}\frac{1}{4}\left(  1+\varepsilon\right)  k^{2}<\frac{1}{1-4\gamma}\left(
\frac{1}{4}+\frac{\gamma}{2}^{2}\right)  k^{2}<\left(  \frac{1}{4}%
+2\gamma\right)  k^{2}.
\]
Thus, we have
\begin{align*}
e\left(  H_{blue}\right)   &  =\binom{k}{2}-e\left(  H_{red}\right)  -e\left(
H_{mid}\right)  -e\left(  H_{irr}\right) \\
&  \geq\binom{k}{2}-\left(  \frac{1}{4}+2\gamma\right)  k^{2}-13\beta
k^{2}-\varepsilon k^{2}>\left(  \frac{1}{4}-16\beta\right)  k^{2}.
\end{align*}
Hence for the size of the graph $G$ we see that
\[
e\left(  G\right)  >e\left(  H_{blue}\right)  \left(  \frac{1}{4}%
-16\beta\right)  k^{2}\left(  1-\beta\right)  t^{2}>\left(  \frac{1}%
{4}-17\beta\right)  n^{2}.
\]
Since we have
\[
bk\left(  G\right)  \leq\left(  \frac{1}{12}-\gamma\right)  n\leq\left(
\frac{1}{6}-2\left(  17\beta\right)  ^{1/3}\right)  n,
\]
by Theorem \ref{stab}, if $\xi$ is sufficiently small then $G$ contains an
induced graph $G_{0}$ with
\[
\left|  G_{0}\right|  \geq\left(  1-\xi\right)  n,
\]
and with
\[
\delta\left(  G_{0}\right)  \geq\left(  1-2\xi\right)  n.
\]
The proof is completed.$\ \hfill\square$

\section{Proof of Theorem \ref{th1}}

Suppose there is some $\varepsilon$ $\left(  0<\varepsilon<1/6\right)  $ such
that for arbitrarily large $q$ there is some $p\leq\left(  1/6-\varepsilon
\right)  q$ such that%
\[
r\left(  B_{p},B_{q}\right)  >2q+3,
\]
i.e. there is a graph $G$ of order $n=2q+3$ such that
\begin{align}
bk\left(  G\right)   &  \leq\left(  \frac{1}{12}-\varepsilon\right)
n,\label{assum1}\\
bk\left(  \overline{G}\right)   &  \leq\frac{n}{2}-2. \label{assum2}%
\end{align}
From Theorem \ref{th2} we see that for every $\xi>0$ if $q$ is large enough
then $G$ must contain an induced bipartite graph $G_{0}$ with
\[
v\left(  G_{0}\right)  \geq\left(  1-\xi\right)  n
\]
and%
\[
\delta\left(  G_{0}\right)  \geq\left(  \frac{1}{2}-2\xi\right)  n.
\]

Let $U=V\left(  G_{0}\right)  $ and $U_{1}$ and $U_{2}$ be the two parts of
$G_{0},$ i.e. $U=U_{1}\cup U_{2}.$ Set
\begin{align*}
V_{0}  &  =V\left(  G\right)  \backslash V\left(  G_{0}\right)  ;\\
V_{1}  &  =\left\{  u:\text{ \ }u\in V_{0},\text{ \ }\Gamma\left(  u\right)
\cap U_{1}\neq\varnothing,\text{ }\Gamma\left(  u\right)  \cap U_{2}%
=\varnothing\right\} \\
V_{2}  &  =\left\{  u:\text{ \ }u\in V_{0},\text{ \ }\Gamma\left(  u\right)
\cap U_{1}=\varnothing,\text{ }\Gamma\left(  u\right)  \cap U_{2}%
\neq\varnothing\right\} \\
V_{3}  &  =\left\{  u:\text{ \ }u\in V_{0},\text{ \ }\Gamma\left(  u\right)
\cap U_{1}\neq\varnothing,\text{ }\Gamma\left(  u\right)  \cap U_{2}%
\neq\varnothing\right\}
\end{align*}
We see immediately that
\begin{align*}
V_{0}  &  =V_{1}\cup V_{2}\cup V_{3}\\
V_{1}\cap V_{2}  &  =\varnothing,\\
e\left(  M_{1},V_{2}\right)   &  =e\left(  M_{2},V_{1}\right)  =0.
\end{align*}
Our goal is to prove that $V_{3}$ is empty.

For every vertex $v\in G_{0}$ set $P\left(  v\right)  =V_{3}\backslash
\Gamma_{G}\left(  v\right)  $ and consider a pair of distinct vertices $u,v\in
U_{2}$. We have
\begin{align*}
bs\left(  \overline{G}\right)   &  \geq\left|  U_{2}\right|  -2+\left|
V_{1}\right|  +\left|  P\left(  u\right)  \cap P\left(  v\right)  \right| \\
&  \geq\left|  U_{2}\right|  -2+\left|  V_{1}\right|  +\left|  P\left(
u\right)  \right|  +\left|  P\left(  v\right)  \right|  -\left|  V_{3}\right|
.
\end{align*}

Summing over all pairs $u,v\in U_{2}$ and taking the average, we obtain%
\begin{align}
bs\left(  \overline{G}\right)   &  \geq\left|  U_{2}\right|  -2+\left|
V_{1}\right|  +\frac{2\left(  \left|  V_{3}\right|  \left|  U_{2}\right|
-e\left(  U_{2},V_{3}\right)  \right)  }{\left|  U_{2}\right|  }-\left|
V_{3}\right| \nonumber\\
&  =\left|  U_{2}\right|  -2+\left|  V_{1}\right|  +\left|  V_{3}\right|
-\frac{2e\left(  U_{2},V_{3}\right)  }{\left|  U_{2}\right|  }. \label{in5}%
\end{align}
As $U_{1}$ is an independent set we have
\[
\left|  U_{1}\right|  \leq bs\left(  \overline{G}\right)  +2\leq\frac{n}{2}.
\]

Thus we obtain,%
\[
\left|  U_{2}\right|  \geq\left|  G_{0}\right|  -\frac{n}{2}>\left(  \frac
{1}{2}-\xi\right)  n.
\]
and hence, estimating $\left|  U_{2}\right|  $ in the denominator in
(\ref{in5}),
\[
bs\left(  \overline{G}\right)  >\left|  U_{2}\right|  -2+\left|  V_{1}\right|
+\left|  V_{3}\right|  -\frac{4e\left(  U_{2},V_{3}\right)  }{\left(
1-2\xi\right)  n}.
\]

By symmetry, we find that
\[
bs\left(  \overline{G}\right)  >\left|  U_{1}\right|  -2+\left|  V_{2}\right|
+\left|  V_{3}\right|  -\frac{4e\left(  U_{1},V_{3}\right)  }{\left(
1-2\xi\right)  n},
\]
and thus, in view of $\left|  U\right|  =\left|  U_{1}\right|  +\left|
U_{2}\right|  $ and $\left|  V_{1}\right|  +\left|  V_{2}\right|  +\left|
V_{3}\right|  +\left|  U\right|  =n,$ yields%
\[
2bs\left(  \overline{G}\right)  >n-4+\left|  V_{3}\right|  -\frac{4e\left(
U,V_{3}\right)  }{\left(  1-2\xi\right)  n}.
\]

By (\ref{assum2}) we immediately obtain
\begin{equation}
e\left(  U,V_{3}\right)  >\frac{\left(  1-2\xi\right)  \left|  V_{3}\right|
n}{4}. \label{in4}%
\end{equation}

On the other hand, every $v\in V_{3}$ has a neighbor $u\in U_{1},$ thus, in
view of
\[
\left|  \Gamma\left(  u\right)  \cap U_{2}\right|  \geq\delta\left(
G_{0}\right)  \geq\left(  \frac{1}{2}-2\xi\right)  n,
\]
we find that%
\begin{align*}
bs\left(  G\right)   &  \geq\left|  \Gamma\left(  v\right)  \cap\Gamma\left(
u\right)  \cap U_{2}\right|  =\left|  \Gamma\left(  v\right)  \cap
U_{2}\right|  +\left|  \Gamma\left(  u\right)  \cap U_{2}\right|  -\left|
U_{2}\right| \\
&  \geq\left|  \Gamma\left(  v\right)  \cap U_{2}\right|  +\left(  \frac{1}%
{2}-2\xi\right)  n-\left|  U_{2}\right|  .
\end{align*}

Taking the average over all $v\in V_{3}$ we obtain%
\[
bs\left(  G\right)  \geq\frac{e\left(  V_{3},U_{2}\right)  }{\left|
V_{3}\right|  }+\left(  \frac{1}{2}-2\xi\right)  n-\left|  U_{2}\right|  ,
\]
and by symmetry,%
\[
bs\left(  G\right)  \geq\frac{e\left(  V_{3},U_{1}\right)  }{\left|
V_{3}\right|  }+\left(  \frac{1}{2}-2\xi\right)  n-\left|  U_{1}\right|  .
\]

Therefore,%
\begin{align*}
2bs\left(  G\right)   &  \geq\frac{e\left(  V_{3},U\right)  }{\left|
V_{3}\right|  }+\left(  1-4\xi\right)  n-\left|  U\right|  \geq\frac{e\left(
V_{3},U\right)  }{\left|  V_{3}\right|  }+\left(  1-4\xi\right)  n-\left|
U\right| \\
&  \geq\frac{e\left(  V_{3},U\right)  }{\left|  V_{3}\right|  }-4\xi n
\end{align*}
and, in view of (\ref{assum1}) we find that,%
\[
\left(  \frac{1}{6}+4\xi\right)  n\left|  V_{3}\right|  \geq e\left(
V_{3},M\right)  .
\]
Combining with (\ref{in4}) we obtain%
\[
\frac{1}{6}+4\xi>\frac{1-2\xi}{4},
\]
and this clearly is a contradiction for small $\xi$. Therefore, $V_{3}%
=\varnothing.$

From $\left|  V_{3}\right|  =0$ we conclude that%
\[
\left|  U_{1}\right|  +\left|  U_{2}\right|  +\left|  V_{1}\right|  +\left|
V_{2}\right|  =n.
\]
To finish the proof it suffices to note that
\begin{align*}
bs\left(  \overline{G}\right)   &  \geq\left|  U_{1}\right|  +\left|
V_{2}\right|  -2,\\
bs\left(  \overline{G}\right)   &  \geq\left|  U_{2}\right|  +\left|
V_{1}\right|  -2,
\end{align*}
and therefore, $bs\left(  \overline{G}\right)  \geq n/2-2$.$\hfill\square$

\section{Tigthness of Theorem \ref{th1}}

We shall show that Theorem \ref{th1} is essentially tight. In particular, we
shall show that if $\epsilon$ is fixed and appropriately small, then for all
sufficiently large $n$ there exists a (partially random) red-bluet coloring of
the edges of $K_{n}$ for which
\[
bk_{R}<\left(  \frac{1}{2}-5\epsilon\right)  n\qquad\text{and}\qquad
bk_{B}<\left(  \frac{1}{12}+3\epsilon\right)  n.
\]

For convenience, assume that $n$ is divisible by $3$. Partition $[n]$ into
three sets $A_{1},A_{2},A_{3}$, each with $n/3$ vertices, and color the graphs
induced by $A_{1},\,A_{2},\,A_{3}$ in red. Then edges of the form $uv$ where
$u\in A_{i}$ and $v\in A_{j}$ $\left(  i\neq j\right)  $ are independently
colored red with probability $p=\frac{1}{2}-\delta$ and blue with probability
$q=\frac{1}{2}+\delta$. For $u,v\in A_{i}$, the size of the red book on $uv$
is a random variable with expected value
\[
\frac{n}{3}-2+\frac{2n}{3}\,p^{2}<\left(  \frac{1}{2}-\frac{2(\delta
-\delta^{2})}{3}\right)  n.
\]

Specifically, the book size is $n/3-2+X$ where $X$ is a Bernoulli random
variable $B(2n/3,p^{2})$. Now suppose $u\in A_{i}$ and $v\in A_{j}$ where
$i\neq j$. If $uv$ is a blue edge, the size of the blue book on $uv$ is a
random variable with expected value
\[
\frac{n}{3}\,q^{2}=\frac{n}{3}\left(  \frac{1}{2}+\delta\right)  ^{2}=\left(
\frac{1}{12}+\frac{\delta+\delta^{2}}{3}\right)  n.
\]

The book size is the Bernoulli random variable $B(n/3,q^{2})$. If $uv$ is a
red edge, the size of the red book on $uv$ is a random variable with expected
value
\[
\frac{n}{3}\,p^{2}+\left(  \frac{2n}{3}-2\right)  p<\left(  \frac{5}%
{12}-\delta+\frac{\delta^{2}}{3}\right)  n.
\]

Set $\delta=8.25\epsilon$. Then
\[
k_{1}=\frac{2(\delta-\delta^{2})}{3}-5\epsilon\text{ \ \ \ \ and
\ \ \ \ }k_{2}=3\epsilon-\frac{(\delta+\delta^{2})}{3}%
\]
are fixed positive numbers. We use the following version of the Chernoff bound
\cite[p. 12]{Bo:RG}: if $X$ is $B(n,p)$ then
\[
P(|X-np|\geq kn)\leq2\exp(-2k^{2}n).
\]

It follows that the probability that that there is a red book with at least
$(\frac{1}{2}-5\epsilon)n$ pages or a blue book with at least $(\frac{1}%
{12}+3\epsilon)n$ pages tends to 0 as $n\rightarrow\infty$. Thus for large
enough $n$ the desired two-coloring of of the edges of $K_{n}$ exists.

This result is easily translated into Ramsey number terms, where it yields the
following fact.

\begin{proposition}
Let $\delta$ be an appropriately small positive number. Then if $q$ is
sufficiently large and $p\geq(\frac{1}{6}+\delta)q$ then
\[
r(B_{p},B_{q})>2(1+\delta)q.
\]
In particular, the relation
\[
r(B_{p},B_{q})=2q+3
\]
no longer holds for $p>(\frac{1}{6}+\delta)q$.
\end{proposition}

\section{\label{secProofs}Proofs}

\subsection{Proof of Lemma \ref{dle}}

We shall prove first that for every $i\in\left[  k\right]  $ there are at
least%
\begin{equation}
d_{i}^{2}\left(  e\left(  A\right)  -2\varepsilon t^{2}\right)  t-2\varepsilon
e\left(  A\right)  t \label{lotri1}%
\end{equation}
triangles in $H$ having $2$ vertices in $A$ and one vertex in $B_{i}.$

This is certianly true if $d_{i}\leq\varepsilon$ as then the above quantity is nonpositive.

Assume $d_{i}>\varepsilon$; applying Lemma \ref{XPle} to the pair $\left(
A,B_{i}\right)  $ we see that there are at most $2\varepsilon t^{2}$ edges
$\left(  u,v\right)  $ in $A$ with%
\[
\left|  \Gamma\left(  u\right)  \cap\Gamma\left(  v\right)  \cap B_{i}\right|
\leq\left(  d_{i}-\varepsilon\right)  ^{2}t.
\]
Hence, there are at least $\left(  e\left(  A\right)  -2\varepsilon
t^{2}\right)  $ edges in $A$ with
\[
\left|  \Gamma\left(  u\right)  \cap\Gamma\left(  v\right)  \cap B_{i}\right|
>\left(  d_{i}-\varepsilon\right)  ^{2}t
\]
and therefore, there are at least
\begin{align*}
\left(  d_{i}-\varepsilon\right)  ^{2}\left(  e\left(  A\right)  -2\varepsilon
t^{2}\right)  t  &  \geq\left(  d_{i}^{2}-2\varepsilon\right)  \left(
e\left(  A\right)  -2\varepsilon t^{2}\right)  t\\
&  >d_{i}^{2}\left(  e\left(  A\right)  -2\varepsilon t^{2}\right)
t-2\varepsilon e\left(  A\right)  t
\end{align*}
triangles in $H$ having $2$ vertices in $A$ and one vertex in $B_{i}.$

Now, summing this inequality for $i=1,...,k$ we obtain the desired
result.$\hfill\square$

\subsection{Proof of Lemma \ref{XPle2}}

Our proof is a straightforward exercise on $\varepsilon$-uniform pairs. Let
$A_{1}^{\prime}$ be the set of all $u\in A_{1}$ such that
\[
\left|  \Gamma\left(  u\right)  \cap B\right|  \leq\left(  d_{1}%
-\varepsilon\right)  \left|  B\right|  .
\]
By the $\varepsilon$-uniformity of $\left(  A_{1},B\right)  $ we have
\[
\left|  A_{1}^{\prime}\right|  \leq\varepsilon\left|  A_{1}\right|
\]
Take any $u_{0}\in A_{1}\backslash A_{1}^{\prime}$ and let $A_{2}^{\prime}$ be
the set of all $v\in A_{2}$ such that
\[
\left|  \Gamma\left(  v\right)  \cap\left(  \Gamma\left(  u\right)  \cap
B\right)  \right|  \leq\left(  d_{2}-\varepsilon\right)  \left|  \Gamma\left(
u\right)  \cap B\right|  .
\]
By the $\varepsilon$-uniformity of $\left(  A_{2},B\right)  ,$ and from
\[
\left|  \Gamma\left(  u\right)  \cap B\right|  >\left(  d_{1}-\varepsilon
\right)  \left|  B\right|  >\varepsilon\left|  B\right|  ,
\]
we obtain%
\[
\left|  A_{2}^{\prime}\right|  \leq\varepsilon\left|  A_{2}\right|  .
\]
Thus, there are at least $\left(  1-\varepsilon\right)  \left|  A_{2}\right|
$ vertices $v\in A_{2}$ such that%
\[
\left|  \Gamma\left(  v\right)  \cap\Gamma\left(  u\right)  \cap B\right|
>\left(  d_{2}-\varepsilon\right)  \left|  \left(  \Gamma\left(  u\right)
\cap B\right)  \right|  >\left(  d_{2}-\varepsilon\right)  \left(
d_{1}-\varepsilon\right)  \left|  B\right|  .
\]
Hence, there are at most
\[
\left(  1-\left(  1-\varepsilon\right)  ^{2}\right)  \left|  A_{1}\right|
\left|  A_{2}\right|  <2\varepsilon\left|  A_{1}\right|  \left|  A_{2}\right|
\]
pairs $\left(  u,v\right)  $ with $u\in A_{1},$ $v\in A_{2}$ such that
(\ref{i03}) holds.$\hfill\square$

\subsection{Proof of Lemma \ref{lebs}}

First we shall prove that for every $i\in\left[  k\right]  $ there are at
least
\begin{equation}
d_{1i}d_{2i}\left(  e\left(  A_{1},A_{2}\right)  -2\varepsilon t^{2}\right)
t-2\varepsilon e\left(  A_{1},A_{2}\right)  t \label{lotri2}%
\end{equation}
triangles in $H$ having one vertex in $A_{1},$ one vertex in $A_{2}$ and one
vertex in $B_{i}.$ This is certianly true if $d_{1i}\leq\varepsilon$ or
$d_{1i}\leq\varepsilon$ as then the above quantity is nonpositive.

Assume $d_{1i}>\varepsilon$ and $d_{1i}>\varepsilon;$ apply Lemma \ref{XPle2}
to the pairs $\left(  A_{1},B_{i}\right)  $ and $\left(  A_{2},B_{i}\right)
.$ Since there are at most $2\varepsilon t^{2}$ pairs $\left(  u,v\right)  $
with $u\in A_{1},$ $v\in A_{2}$ with%
\[
\left|  \Gamma\left(  u\right)  \cap\Gamma\left(  v\right)  \cap B_{i}\right|
\leq\left(  d_{1i}-\varepsilon\right)  \left(  d_{2i}-\varepsilon\right)  t,
\]
there are at least $\left(  e\left(  A_{1},A_{2}\right)  -2\varepsilon
t^{2}\right)  $ edges $\left(  u,v\right)  \in E\left(  A_{1},A_{2}\right)  $
with
\[
\left|  \Gamma\left(  u\right)  \cap\Gamma\left(  v\right)  \cap B_{i}\right|
>\left(  d_{1i}-\varepsilon\right)  \left(  d_{2i}-\varepsilon\right)  t.
\]
Therefore, there are at least
\begin{align*}
\left(  d_{1i}-\varepsilon\right)  \left(  d_{2i}-\varepsilon\right)  \left(
e\left(  A_{1},A_{2}\right)  -2\varepsilon t^{2}\right)  t  &  \geq\left(
d_{1i}d_{2i}-2\varepsilon\right)  \left(  e\left(  A_{1},A_{2}\right)
-2\varepsilon t^{2}\right)  t\\
&  \geq d_{1i}d_{2i}\left(  e\left(  A_{1},A_{2}\right)  -2\varepsilon
t^{2}\right)  t\\
&  -2\varepsilon e\left(  A_{1},A_{2}\right)  t
\end{align*}
triangles in $H$ having one vertex in $A_{1},$ one vertex in $A_{2}$ and one
vertex in $B_{i}.$

Now, summing this inequality for all $i\in\left[  k\right]  $ we obtain the
desired result.$\hfill\square$

\subsection{Proof of Claim \ref{cl1}}

Assume the opposite and let $V_{i}$ be a blob with%
\[
e_{R}\left(  V_{i}\right)  \geq\gamma t^{2},\text{ \ \ \ \ }e_{B}\left(
V_{i}\right)  \geq\gamma t^{2}%
\]
Let us compute the average size of the blue books whose base is $E\left(
V_{i}\right)  .$ Let $M\subset\left[  k\right]  $ be the set of all
$s\in\left[  k\right]  \backslash\left\{  i\right\}  $ such that $\left(
V_{i},V_{s}\right)  $ is an $\varepsilon$-uniform pair; by our assumption
$\left|  M\right|  \geq\left(  1-\varepsilon\right)  k$. Applying Corollary
\ref{bs1} with
\[
A=V_{i},\text{ }B_{s}=V_{s}:s\in M
\]
we obtain
\[
bk_{B}\geq t\left(  1-\gamma\right)  \sum_{s\in M}\left(  1-d_{is}\right)
^{2}-2\varepsilon kt,\text{ }%
\]
and by Cauchy's inequality, in view of $\left|  M\right|  \leq k$ and
(\ref{assum}), we see that%
\begin{align*}
\frac{1}{k}\left(  1-\gamma\right)  \left(  \sum_{s\in M}\left(
1-d_{is}\right)  \right)  ^{2}  &  \leq\left(  \frac{1}{12}-\gamma
+2\varepsilon\right)  \left(  1+\varepsilon\right)  k\\
&  \leq\left(  \frac{1}{12}-\gamma+\gamma^{2}\right)  \left(  1+\gamma
^{2}/2\right)  k
\end{align*}
Hence,
\begin{equation}
\sum_{s\in M}\left(  1-d_{is}\right)  \leq\sqrt{\left(  \frac{1}{12}%
-\gamma+\gamma^{2}\right)  \frac{\left(  1+\gamma^{2}/2\right)  }{\left(
1-\gamma\right)  }}k\leq\sqrt{\left(  \frac{1}{12}-\frac{\gamma}{2}\right)
}k. \label{bbmax}%
\end{equation}
Similarly, estimating the average size of the red books whose base is in
$V_{i},$ we see that%
\[
\frac{1}{k}\left(  \sum_{s\in M}d_{is}\right)  ^{2}\leq\left(  \frac{1}%
{2}+2\varepsilon\right)  \frac{1+\varepsilon}{1-\gamma}k=\left(  \frac{1}%
{2}+\gamma^{2}\right)  \frac{\left(  1+\gamma^{2}/2\right)  }{\left(
1-\gamma\right)  }k\leq\left(  \frac{1}{2}+\gamma\right)  k
\]
and hence,
\begin{equation}
\sum_{s\in M}d_{is}\leq k\sqrt{\frac{1}{2}+\gamma} \label{rbmax}%
\end{equation}
Adding (\ref{bbmax}) and (\ref{rbmax}) we find that
\begin{align*}
\sqrt{\left(  \frac{1}{2}+\gamma\right)  }+\sqrt{\left(  \frac{1}{12}%
-\frac{\gamma}{2}\right)  }  &  \geq\frac{\left|  M\right|  }{k}\\
&  \geq1-\varepsilon>1-\gamma.
\end{align*}
This leads to a contradiction if $\gamma$ is small, as
\[
\frac{1}{\sqrt{2}}+\sqrt{\frac{1}{12}}=\frac{3\sqrt{2}+\sqrt{3}}%
{6}=0.995...<1.
\]

$\hfill\square$

\subsection{Proof of Claim \ref{cl2}}

Assume three blobs satisfying the conditions of the claim exist. Let $M$ be
the set of all $s\in\left[  k\right]  \backslash\left\{  i,j,l\right\}  $ such
that every one of the pairs $\left(  V_{i},V_{s}\right)  ,$ $\left(
V_{j},V_{s}\right)  ,$ $\left(  V_{l},V_{s}\right)  $ is $\varepsilon
$-uniform; clearly $\left|  M\right|  \geq\left(  1-3\varepsilon\right)  k$.

As in the proof of the Claim \ref{cl1}, estimating the average size of the red
books with base in one of the sets $E\left(  V_{i}\right)  ,E\left(
V_{j}\right)  ,E\left(  V_{l}\right)  $ gives
\begin{align}
bk_{R}  &  \geq t\left(  1-\gamma\right)  \sum_{s\in M}d_{is}^{2}%
-2kt\varepsilon,\label{bs1.1}\\
bk_{R}  &  \geq t\left(  1-\gamma\right)  \sum_{s\in M}d_{js}^{2}%
-2kt\varepsilon,\label{bs1.2}\\
bk_{R}  &  \geq t\left(  1-\gamma\right)  \sum_{s\in M}d_{ls}^{2}%
-2kt\varepsilon. \label{bs1.3}%
\end{align}
On the other hand, applying Corollary \ref{bs2} with
\[
A_{1}=V_{i},\text{ }A_{2}=V_{j},\text{ }B_{s}=V_{s}:s\in M
\]
we obtain for the average size $S$ of the blue books with base in $E\left(
V_{i},V_{j}\right)  $
\begin{align}
bk_{B}  &  \geq S\geq t\left(  1-\frac{2\varepsilon t^{2}}{e_{B}\left(
A_{1},A_{2}\right)  }\right)  \sum_{s\in M}\left(  1-d_{is}\right)  \left(
1-d_{js}\right)  -2\varepsilon kt\nonumber\\
&  \geq t\left(  1-\gamma\right)  \sum_{s\in M}\left(  1-d_{is}\right)
\left(  1-d_{js}\right)  -2\varepsilon kt \label{bs2.1}%
\end{align}
Considering in turn $\left(  V_{i},V_{k}\right)  $ and $\left(  V_{j}%
,V_{k}\right)  $ we obtain exactly in the same way%
\begin{align}
bk_{B}  &  \geq t\left(  1-\gamma\right)  \sum_{s\in M}\left(  1-d_{is}%
\right)  \left(  1-d_{ls}\right)  -2\varepsilon kt\label{bs2.2}\\
bk_{B}  &  \geq t\left(  1-\gamma\right)  \sum_{s\in M}\left(  1-d_{js}%
\right)  \left(  1-d_{ls}\right)  -2\varepsilon kt. \label{bs2.3}%
\end{align}
Setting
\[
d_{s}=\sum_{s\in M}d_{is}+d_{js}+d_{ls}%
\]
and adding (\ref{bs2.1}), (\ref{bs2.2}), (\ref{bs2.3}) together with each of
(\ref{bs1.1}), (\ref{bs1.2}), (\ref{bs1.3}) multiplied by 1/2, we obtain
\[
t\left(  1-\gamma\right)  \sum_{s\in M}\left(  3-2d_{s}+\frac{1}{2}d_{s}%
^{2}\right)  -9\varepsilon kt\leq3bk_{B}+\frac{3}{2}bk_{R}.
\]
Hence, setting
\[
d=\frac{1}{\left|  M\right|  }\sum_{s\in M}d_{s},
\]
by Cauchy's inequality and (\ref{assum}), we see that
\[
\left|  M\right|  t\left(  1-\gamma\right)  \left(  3-2d+\frac{1}{2}%
d^{2}\right)  -9\varepsilon kt\leq\left(  3\left(  \frac{1}{12}-\gamma\right)
+\frac{3}{4}\right)  \left(  1+\varepsilon\right)  kt.
\]
Therefore,
\[
\left(  1-\gamma\right)  \left(  1-3\varepsilon\right)  \left(  3-2d+\frac
{1}{2}d^{2}\right)  -9\varepsilon\leq\left(  1-3\gamma\right)  \left(
1+\varepsilon\right)  .
\]
Hence, by $\varepsilon=\gamma^{2}/2,$
\begin{align*}
2-2d+\frac{1}{2}d^{2}  &  <\frac{\left(  1-3\gamma\right)  \left(
1+\gamma^{2}/2\right)  }{\left(  1-\gamma\right)  \left(  1-3\gamma
^{2}/2\right)  }+8\gamma^{2}-1\\
&  =\left(  1-\frac{3\gamma}{\left(  1-\gamma\right)  }\right)  \left(
1+\frac{2\gamma^{2}}{1-3\gamma^{2}/2}\right)  +8\gamma^{2}-1\\
&  <-\frac{3\gamma}{\left(  1-\gamma\right)  }+10\gamma^{2}.
\end{align*}
This is a contradiction since the right-hand side is negative for $\gamma<1$
while the left-hand side is always nonnegative.$\hfill\square$

\subsection{Proof of Claim \ref{cl3}}

Assume two blobs satisfying the conditions of the claim exist. Let $M$ be the
set of all $s\in\left[  k\right]  \backslash\left\{  i,j\right\}  $ such that
every one of the pairs $\left(  V_{i},V_{s}\right)  ,\left(  V_{j}%
,V_{s}\right)  $ is $\varepsilon$-uniform; clearly $\left|  M\right|
\geq\left(  1-2\varepsilon\right)  k$.

As in the proof of the Claim \ref{cl1} estimating the average size of the blue
books with base in one of the sets $E\left(  V_{i}\right)  $ or $E\left(
V_{j}\right)  $ we obtain
\begin{align}
bk_{B}  &  \geq t\left(  1-\gamma\right)  \sum_{s\in M}\left(  1-d_{is}%
\right)  ^{2}-2\varepsilon kt,\label{bs3.1}\\
bk_{B}  &  \geq t\left(  1-\gamma\right)  \sum_{s\in M}\left(  1-d_{js}%
\right)  ^{2}-2\varepsilon kt. \label{bs3.2}%
\end{align}

As in the proof of the Claim \ref{cl2} by estimating the average size of the
red books having a base in $E\left(  V_{i},V_{j}\right)  $ we obtain%
\begin{equation}
bk_{R}\geq t\left(  1-\gamma\right)  \sum_{s\in M}d_{is}d_{js}-2\varepsilon
kt. \label{bs4}%
\end{equation}

Setting
\[
d_{s}=\sum_{s\in M}d_{is}+d_{js},
\]
and adding (\ref{bs3.1}),(\ref{bs3.2}), and doubled (\ref{bs4}) we obtain%
\[
t\left(  1-\gamma\right)  \sum_{s\in M}\left(  2-2d_{s}+d_{s}^{2}\right)
-8\varepsilon kt\leq2bk_{B}+2bk_{R}.
\]
Hence, by (\ref{assum}) we see that
\begin{align*}
\left(  1-\gamma\right)  \sum_{s\in M}\left(  2-2d_{s}+d_{s}^{2}\right)   &
\leq\left(  2\left(  \frac{1}{12}-\gamma\right)  +1\right)  \left(
1+\varepsilon\right)  k+8\varepsilon k\\
&  =\left(  \frac{7}{6}-2\gamma+\frac{7}{6}\varepsilon-2\gamma\varepsilon
+8\varepsilon\right)  k<\frac{7}{6}k.
\end{align*}
Setting
\[
d=\frac{1}{\left|  M\right|  }\sum_{s\in M}d_{s},
\]
by Cauchy's inequality and $\left|  M\right|  \geq\left(  1-2\varepsilon
\right)  k,$
\[
\left(  1-2\varepsilon\right)  k\left(  2-2d+d^{2}\right)  \leq\frac{7}%
{6}\frac{1}{1-\gamma}k.
\]
and hence%
\begin{equation}
2-2d+d^{2}\leq\frac{7}{6}\frac{1}{\left(  1-2\varepsilon\right)  \left(
1-\gamma\right)  }. \label{din1}%
\end{equation}
On the other hand, applying Cauchy's inequality to (\ref{bs3.1}) and
(\ref{bs3.2}), in view of (\ref{assum}), we obtain%
\begin{align*}
\left|  M\right|  \left(  1-\frac{1}{\left|  M\right|  }\sum_{s\in M}%
d_{is}\right)  ^{2}-2\varepsilon k  &  \leq\left(  \frac{1}{12}-\gamma\right)
\left(  1+\varepsilon\right)  k,\text{\ }\\
\left|  M\right|  \left(  1-\frac{1}{\left|  M\right|  }\sum_{s\in M}%
d_{js}\right)  ^{2}-2\varepsilon k  &  \leq\left(  \frac{1}{12}-\gamma\right)
\left(  1+\varepsilon\right)  k.
\end{align*}
Hence again by Cauchy's inequality, and $\left|  M\right|  \geq\left(
1-2\varepsilon\right)  k,$ we see that
\[
\left(  1-2\varepsilon\right)  \left(  1-\frac{1}{2}d\right)  ^{2}\leq\left(
\frac{1}{12}-\gamma\right)  \left(  1+\varepsilon\right)  +4\varepsilon,
\]
yielding%
\begin{equation}
\left(  2-d\right)  ^{2}\leq\left(  \frac{1}{3}-\gamma\right)  \frac
{1+\gamma^{2}/2}{1-\gamma^{2}}+\frac{2\gamma^{2}}{1-\gamma^{2}}. \label{din2}%
\end{equation}
Since we can select $\gamma$ arbitrarily small, from (\ref{din1}) and
(\ref{din2}) we obtain
\[
2-\frac{1}{\sqrt{3}}\leq d\leq1+\sqrt{\frac{1}{6}}%
\]
giving%
\[
1\leq\frac{2\sqrt{3}+\sqrt{6}}{6},
\]
a contradiction, since
\[
\frac{2\sqrt{3}+\sqrt{6}}{6}=0.98...<1.
\]

$\hfill\square$

\subsection{Proof of Claim \ref{cl4}}

The proof of this claim is by far the most complicated one.

Assume $V_{1},...,V_{l}$ are the blue blobs and set $l=\alpha k.$

Fix a blue blob $V_{i}$ and estimate the average size of the blue books whose
base is $E\left(  V_{i}\right)  $. Let $M_{i}$ be the set of all $s\in\left[
k\right]  \backslash\left\{  i\right\}  $ such that $\left(  V_{i}%
,V_{s}\right)  $ is an $\varepsilon$-uniform pair and let
\begin{align*}
M_{i1}  &  =M_{i}\cap\left[  l\right]  ,\\
\text{ }M_{i2}  &  =M_{i}\cap\left[  l+1,k\right]  .
\end{align*}
From $\left|  M_{i}\right|  \geq\left(  1-\varepsilon\right)  k$ we
immediately obtain%
\begin{align}
l  &  \geq\left|  M_{i1}\right|  \geq l-\varepsilon k,\label{m1}\\
k-l  &  \geq\left|  M_{i2}\right|  \leq k-l-\varepsilon k.\nonumber
\end{align}
By Claim \ref{cl3} for every $s\in M_{i1}$ the red density of the pair
$\left(  V_{i},V_{s}\right)  $ satisfies%
\[
1-d_{is}\geq1-\gamma.
\]
Therefore, the average size $S$ of the blue books whose base is in $V_{i}$
satisfies%
\begin{align*}
S  &  \geq\left(  \left(  1-\varepsilon-\gamma\right)  ^{2}\left|
M_{i1}\right|  +\sum_{s\in M_{i2}}\left(  1-d_{is}\right)  ^{2}-2\varepsilon
k\right)  t\\
&  \geq\left(  \left(  1-4\gamma\right)  \left|  M_{i1}\right|  +\sum_{s\in
M_{i2}}\left(  1-d_{is}\right)  ^{2}-2\varepsilon k\right)  t.
\end{align*}

Hence, by $bk_{B}\geq S,$ in view of (\ref{assum}), we see that
\[
\left(  1-4\gamma\right)  \left|  M_{i1}\right|  +\sum_{s\in M_{i2}}\left(
1-d_{is}\right)  ^{2}\leq\left(  \frac{1}{12}-\gamma+2\varepsilon\right)
\left(  1+\varepsilon\right)  k.
\]

Then, by Cauchy's inequality and (\ref{m1}), we find that%
\[
\left(  1-4\gamma\right)  \left(  l-\varepsilon k\right)  +\frac{1}%
{k-l}\left(  \sum_{s\in M_{i2}}1-d_{is}\right)  ^{2}\leq\left(  \frac{1}%
{12}-\gamma+2\varepsilon\right)  \left(  1+\varepsilon\right)  k,
\]
and hence
\begin{align*}
\sum_{s\in M_{i2}}1-d_{is}  &  \leq k\sqrt{\left(  \left(  \frac{1}{12}%
-\gamma+2\varepsilon\right)  \left(  1+\varepsilon\right)  -\left(
1-4\gamma\right)  \left(  \alpha-\varepsilon\right)  \right)  \left(
1-\alpha\right)  }\\
&  .
\end{align*}
Summing this inequality for $i=1,...,l$ and setting%
\[
A=\alpha\sqrt{\left(  \left(  \frac{1}{12}-\gamma+2\varepsilon\right)  \left(
1+\varepsilon\right)  -\left(  1-4\gamma\right)  \left(  \alpha-\varepsilon
\right)  \right)  \left(  1-\alpha\right)  }%
\]
we obtain
\begin{equation}
\sum_{i\in\left[  l\right]  }\sum_{s\in M_{i2}}\left(  1-d_{is}\right)  \leq
Ak^{2}. \label{bbl}%
\end{equation}

Next we shall obtain a similar inequality by considering the average size of
the red books whose bases are contained in a red blob.

Let us define the graph $H$ as follows. The vertices of $H$ are the numbers
$\left[  l+1,k\right]  $ and two vertices $i,j$ are joined iff the red density
$d_{ij}$ of the pair $\left(  V_{i},V_{j}\right)  $ satisfies
\[
d_{ij}\geq1-\gamma.
\]
By Claim \ref{cl2} the complement of $H$ is triangle-free, hence, by
Tur\'{a}n's theorem, the size of $H$ satisfies
\begin{equation}
e\left(  H\right)  \geq\left(  \frac{1}{4}-\varepsilon\right)  \left(
k-l\right)  ^{2} \label{hs}%
\end{equation}
if $k$ is sufficiently large. Fix some $i\in\left[  l+1,k\right]  .$ Let
$M_{i}$ be the set of all $s\in\left[  k\right]  \backslash\left\{  i\right\}
$ such that $\left(  V_{i},V_{s}\right)  $ is an $\varepsilon$-uniform pair
and let
\begin{align*}
M_{i1}  &  =M_{i}\cap\Gamma_{H}\left(  i\right)  ,\\
M_{i2}  &  =M_{i}\cap\left[  l\right]  .
\end{align*}
From $\left|  M_{i}\right|  \geq\left(  1-\varepsilon\right)  k$ we
immediately obtain%
\begin{align}
d_{H}\left(  i\right)   &  \geq\left|  M_{i1}\right|  \geq d_{H}\left(
i\right)  -\varepsilon k,\label{m2}\\
l  &  \geq\left|  M_{i2}\right|  \geq l-\varepsilon k.\nonumber
\end{align}

Therefore, the average size $S$ of the red books whose base is in $E\left(
V_{i}\right)  $ satisfies%
\begin{align*}
S  &  \geq\left(  \left(  1-\varepsilon-\gamma\right)  ^{2}\left|
M_{i1}\right|  +\sum_{s\in M_{i2}}d_{is}^{2}-2\varepsilon k\right)  t\\
&  \geq\left(  \left(  1-4\gamma\right)  \left|  M_{i1}\right|  +\sum_{s\in
M_{i2}}d_{is}^{2}-2\varepsilon k\right)  t.
\end{align*}

Hence, by $bk_{R}\geq S,$ in view of (\ref{assum}) and (\ref{m2}),\ we see
that
\[
\left(  1-4\gamma\right)  \left(  d_{H}\left(  i\right)  -\varepsilon
k\right)  +\sum_{s\in M_{i2}}d_{is}^{2}\leq\frac{1}{2}\left(  1+5\varepsilon
\right)  k.
\]

Thus, by Cauchy's inequality and (\ref{m2}), we find that%
\[
\left(  1-4\gamma\right)  \left(  d_{H}\left(  i\right)  -\varepsilon
k\right)  +\frac{1}{l}\left(  \sum_{s\in M_{i2}}d_{is}\right)  ^{2}\leq
\frac{1}{2}\left(  1+5\varepsilon\right)  k.
\]

Summing this inequality for $i=l+1,...,k$ we obtain%
\begin{align*}
\frac{1}{l}\sum_{i=l+1}^{k}\left(  \sum_{s\in M_{i2}}d_{is}\right)  ^{2}  &
\leq\frac{1}{2}\left(  1+5\varepsilon\right)  k\left(  k-l\right) \\
&  -\left(  1-4\gamma\right)  \left(  2e\left(  H\right)  -\varepsilon
k\left(  k-l\right)  \right)  .
\end{align*}
Hence by (\ref{hs}) we see that%
\begin{align*}
\frac{1}{l}\sum_{i=l+1}^{k}\left(  \sum_{s\in M_{i2}}d_{is}\right)  ^{2}  &
\leq\frac{1}{2}\left(  1+5\varepsilon\right)  k\left(  k-l\right) \\
&  -\left(  1-4\gamma\right)  \left(  \left(  \frac{1}{2}-2\varepsilon\right)
\left(  k-l\right)  -\varepsilon k\right)  \left(  k-l\right)  .
\end{align*}

Applying Cauchy's inequality and replacing $l$ by $\alpha k$ we obtain%
\begin{align*}
\left(  \sum_{i=l+1}^{k}\sum_{s\in M_{i2}}d_{is}\right)  ^{2}  &  \leq\frac
{1}{2}\left(  1+5\varepsilon\right)  \alpha\left(  1-\alpha\right)  ^{2}%
k^{4}\\
&  -\left(  \left(  1-4\gamma\right)  \left(  \left(  \frac{1}{2}%
-2\varepsilon\right)  \left(  1-\alpha\right)  -\varepsilon\right)  \right)
\alpha\left(  1-\alpha\right)  ^{2}k^{4}.
\end{align*}

Setting for brevity%
\[
B=\left(  1-\alpha\right)  \sqrt{\left(  \frac{1}{2}\left(  1+5\varepsilon
\right)  -\left(  \left(  1-4\gamma\right)  \left(  \left(  \frac{1}%
{2}-2\varepsilon\right)  \left(  1-\alpha\right)  -\varepsilon\right)
\right)  \right)  \alpha}%
\]
we obtain%
\begin{equation}
\sum_{i=l+1}^{k}\sum_{s\in M_{i2}}\left(  d_{is}-\varepsilon\right)  \leq
Bk^{2}. \label{bbl1}%
\end{equation}

Now adding (\ref{bbl}) to (\ref{bbl1}) we have%
\[
\left(  1-2\varepsilon\right)  \sum_{i=l+1}^{k}\sum_{s\in M_{i2}}\leq\left(
A+B\right)  k^{2}.
\]

Observe that the sum
\[
\sum_{i=l+1}^{k}\sum_{s\in M_{i2}}%
\]
is just the number of the $\varepsilon$-uniform pairs joining blue to red
blobs, and hence,
\[
\sum_{i=l+1}^{k}\sum_{s\in M_{i2}}\geq l\left(  k-l\right)  -\varepsilon
k^{2}=\left(  \left(  1-\alpha\right)  \alpha-\varepsilon\right)  k^{2}.
\]

Thus, we see that
\[
\left(  1-2\varepsilon\right)  \left(  \left(  1-\alpha\right)  \alpha
-\varepsilon\right)  \leq A+B.
\]

If we assume that this inequality holds for arbitrary small $\gamma$, we
obtain%
\[
\left(  1-\alpha\right)  \alpha\leq\left(  1-\alpha\right)  \alpha\frac
{1}{\sqrt{2}}+\alpha\sqrt{\left(  \frac{1}{12}-\alpha\right)  \left(
1-\alpha\right)  }%
\]
and consequently,%
\[
\sqrt{\left(  1-\alpha\right)  }\left(  1-\frac{1}{\sqrt{2}}\right)  \leq
\sqrt{\left(  \frac{1}{12}-\alpha\right)  },
\]
implying%
\[
\alpha\left(  \sqrt{2}-\frac{1}{2}\right)  \leq\sqrt{2}-\frac{17}%
{12}=-.002...<0,
\]
which is a contradiction.$\hfill\square$

\end{document}